# NONLINEAR SCHRÖDINGER EQUATIONS WITH STARK POTENTIAL

RÉMI CARLES AND YOSHIHISA NAKAMURA

ABSTRACT. We study the nonlinear Schrödinger equations with a linear potential. A change of variables makes it possible to deduce results concerning finite time blow up and scattering theory from the case with no potential.

## 1. INTRODUCTION

In this note, we consider the nonlinear Schrödinger equation with Stark effect,

$$(1.1) \quad \begin{cases} i\varepsilon\partial_t u + \frac{1}{2}\varepsilon^2 \Delta u = V(x)u + \lambda |u|^{2\sigma} u \ , \\ u_{|t=0} = u_0 \ , \end{cases}$$

where $x \in \mathbb{R}^n$, and the potential $V$ is linear,

$$(1.2) \quad V(x) = E \cdot x \ ; \ E = (E_1, \ldots, E_n) \in \mathbb{R}^n \setminus \{0\} \ .$$

We assume that $\varepsilon \in ]0,1]$, $\lambda \in \mathbb{R}$, $\sigma > 0$, and $\sigma < 2/(n-2)$ if $n \geq 3$.

The vector $E$ may represent a constant electric field (see e.g. [6]), or gravity (see e.g. [25]). We introduce the factor $\varepsilon$ to treat both the quantum case where $\varepsilon = \hbar$ (see e.g. [25]), and the case $\varepsilon = 1$, where the nonlinear Schrödinger equation may appear as an envelope equation (see e.g. [24], [8], [13]).

We compare solutions of (1.1) with solutions to the nonlinear Schrödinger equation,

$$(1.3) \quad \begin{cases} i\varepsilon\partial_t v + \frac{1}{2}\varepsilon^2 \Delta v = \lambda |v|^{2\sigma} v \ , \\ v_{|t=0} = u_0 \ . \end{cases}$$

Many results on (1.3) are available, in particular when $u_0 \in L^2(\mathbb{R}^n)$, $H^1(\mathbb{R}^n)$, or

$$\Sigma := \left\{ f \in H^1(\mathbb{R}^n) \ ; \ |x|f \in L^2(\mathbb{R}^n) \right\} \ .$$

Some existence results for (1.1) in the case $\varepsilon = 1$ are consequences of results stated in [7] and [17] (see [18]). Natural questions arise, such as local/global existence issues, finite time blow up and scattering theory. We restrict our study to these directions, but other results could be stated easily, using Proposition 3.1 below. This proposition shows that an explicit change of variables, known as Avron-Herbst formula in the linear setting, makes the comparison between solutions to (1.1) and solutions to (1.3) very easy for the questions mentioned above. On the other hand, Proposition 3.1 does not seem useful to study solitary waves.

2000 *Mathematics Subject Classification.* 35Q55, 35B05, 35P25.

This work was partially supported by the ACI grant "Équation des ondes : oscillations, dispersion et contrôle" (first author) and the JSPS Research Fellowships for Young Scientists (second author).

We first show in Sect. 2 that solutions to (1.1) satisfy a pseudo-conformal evolution law which is the analog to the case $E = 0$, derived in [11]. We then turn to the main result of this note in Sect. 3.

## 2. Some preliminary results

The semiclassical limit for (1.1) was studied in [3] (for more general potentials). It was noticed that the following operator, introduced in [22],

$$J_E^\varepsilon(t) = \frac{x}{\varepsilon} + it\nabla_x - \frac{t^2}{2\varepsilon}E \;,$$

enjoys properties which are agreeable to study the nonlinear equation (1.1). In [22] and [3], this operator is introduced as an Heisenberg observable, therefore the first property in the following lemma is obvious. The other properties were noticed in [3], and they suit nonlinear problems. The last two points follow from the second one.

**Lemma 2.1** ([3], Lemma 2.1). *The operator $J_E^\varepsilon$ is such that:*

- *It commutes with the linear part of (1.1),*

(2.1) $$\left[ J_E^\varepsilon(t), i\varepsilon\partial_t + \frac{1}{2}\varepsilon^2\Delta - V(x) \right] = 0 \;.$$

- *Denote*

$$\phi(t,x) := \frac{|x|^2}{2t} - \frac{t}{2}E\cdot x - \frac{t^3}{24}|E|^2 \;.$$

*Then $\phi$ solves the eikonal equation*

(2.2) $$\partial_t \phi + \frac{1}{2}|\nabla_x \phi|^2 + V(x) = 0 \;,$$

*and for $t \neq 0$,*

(2.3) $$J_E^\varepsilon(t) = ite^{i\phi(t,x)/\varepsilon} \nabla_x \left( e^{-i\phi(t,x)/\varepsilon} \cdot \right).$$

- *For $r \geq 2$, and $r < \frac{2n}{n-2}$ if $n \geq 3$ ($r \leq \infty$ if $n = 1$), define $\delta(r)$ by*

$$\delta(r) \equiv n\left( \frac{1}{2} - \frac{1}{r} \right).$$

*There exists $C_r$ such that, for any $f \in \Sigma$,*

(2.4) $$\|f\|_{L^r} \leq \frac{C_r}{|t|^{\delta(r)}} \|f\|_{L^2}^{1-\delta(r)} \|J_E^\varepsilon(t)f\|_{L^2}^{\delta(r)}.$$

- *For any function $F \in C^1(\mathbb{C}, \mathbb{C})$ satisfying the gauge invariance condition*

$$\exists G \in C(\mathbb{R}_+, \mathbb{R}), \; F(z) = zG(|z|^2),$$

*one has, for $t \neq 0$,*

(2.5) $$J_E^\varepsilon(t)F(w) = \partial_z F(w) J_E^\varepsilon(t) w - \partial_{\bar{z}} F(w) \overline{J_E^\varepsilon(t)w}.$$

Assume that $\varepsilon = 1$. After [1] and [2], we may hope that the evolution of $\|J_E(t)u\|_{L^2}$ is described by a law similar to the pseudo-conformal conservation law, derived in [11]. Indeed, we have the following preliminary result.

**Lemma 2.2.** *Let $u$ solve (1.1) with $\varepsilon = 1$. Then formally, the evolution of $\|J_E(t)u\|_{L^2}$ is given by the pseudo-conformal conservation law,*

$$(2.6) \quad \frac{d}{dt}\left(\frac{1}{2}\|J_E(t)u\|_{L^2}^2 + \frac{\lambda t^2}{\sigma+1}\|u(t)\|_{L^{2\sigma+2}}^{2\sigma+2}\right) = \frac{\lambda t}{\sigma+1}(2-n\sigma)\|u(t)\|_{L^{2\sigma+2}}^{2\sigma+2}.$$

*Proof.* The proof we give is also valid for the case with no potential, and uses (2.1) and (2.5) (recall that this last property is a consequence of (2.3)).

First, let $(e_1, \ldots, e_n)$ be an orthonormal basis of $\mathbb{R}^n$ such that $e_1$ is parallel to $E$, in which we have $V(x) = |E|x_1$. Since the result we want to prove is intrinsic and does not depend on the chosen basis, it is enough to prove it in this context.

Apply the operator $J_E(t)$ to (1.1). Thanks to (2.1), we have

$$i\partial_t J_E(t)u + \frac{1}{2}\Delta J_E(t)u = V(x)J_E(t)u + \lambda J_E(t)\left(|u|^{2\sigma}u\right).$$

The last term can be simplified, thanks to (2.5). Multiply the above equation by $\overline{J_E(t)u}$, integrate with respect to the space variable, and take the imaginary part. This yields,

$$\frac{1}{2}\frac{d}{dt}\|J_E(t)u\|_{L^2}^2 = -\lambda\sigma\,\mathrm{Im}\int u^2|u|^{2\sigma-2}\left(\overline{J_E(t)u}\right)^2 dx.$$

Expanding,

$$(J_E(t)u)^2 = |x|^2 u^2 + \frac{|E|^2 t^4}{4}u^2 - |E|t^2 x_1 u + 2itux\cdot\nabla_x u - t^2(\nabla_x u)^2 - i|E|t^3 u\partial_1 u.$$

The first three terms vanish when taking the imaginary part. The next two do not see the potential, so if the last one vanishes, the lemma follows from the classical pseudo-conformal conservation law. We have

$$\mathrm{Im}\int u^2|u|^{2\sigma-2}i\overline{u}\partial_1\overline{u}dx = \mathrm{Re}\int u|u|^{2\sigma}\partial_1\overline{u}dx = \frac{1}{2}\int |u|^{2\sigma}\partial_1|u|^2 dx = 0.$$

The end of the proof consists in computing the time derivative of $\|u(t)\|_{L^{2\sigma+2}}^{2\sigma+2}$,

$$\frac{d}{dt}\|u(t)\|_{L^{2\sigma+2}}^{2\sigma+2} = (2\sigma+2)\,\mathrm{Re}\int |u|^{2\sigma}\overline{u}\partial_t u\,dx = (\sigma+1)\,\mathrm{Im}\int |u|^{2\sigma}\overline{u}\Delta u\,dx,$$

where we used the fact that the potential is real. The end of the proof is easy computations (the dependence upon the potential has disappeared), we leave out this part. □

## 3. The main results

The fact that we get exactly the same law as in the case with no potential is puzzling. It is explained by the following result.

**Proposition 3.1.** *Let $v^\varepsilon$ solve the initial value problem (1.3). Define $u^\varepsilon$ by*

$$(3.1) \quad u^\varepsilon(t,x) = v^\varepsilon\left(t, x + \frac{t^2}{2}E\right)e^{-i\left(tE\cdot x + \frac{t^3}{6}|E|^2\right)/\varepsilon}.$$

*Then $u^\varepsilon$ solves (1.1). Conversely, if $u^\varepsilon$ solves (1.1), then $v^\varepsilon$, defined by*

$$(3.2) \quad v^\varepsilon(t,x) = u^\varepsilon\left(t, x - \frac{t^2}{2}E\right)e^{i\left(tE\cdot x - \frac{t^3}{3}|E|^2\right)/\varepsilon},$$

*solves (1.3). In particular, if the solution to (1.3) is unique, then so is the solution to (1.1), which is given by (3.1).*

This result was known in the linear case ($\lambda = 0$) (Avron-Herbst formula, see e.g. [6], Chapter 7); we investigate some of its consequences in a nonlinear setting.

*Remark.* Of course, this result would still hold if the nonlinearity $\lambda |z|^{2\sigma} z$ was replaced by a more general one, of the form $G(|z|)z$, since the argument of $G$ is spherically invariant.

*Remark.* The result also holds when an interaction of Hartree type is present. Formula (3.2) turns solutions of

$$i\varepsilon \partial_t u + \frac{1}{2}\varepsilon^2 \Delta u = E \cdot x u + \lambda |u|^{2\sigma} u + \mu \left( |x|^{-\gamma} * |u|^2 \right) u \; ; \; u_{|t=0} = u_0 ,$$

into solutions of

$$i\varepsilon \partial_t v + \frac{1}{2}\varepsilon^2 \Delta v = \lambda |v|^{2\sigma} v + \mu \left( |x|^{-\gamma} * |v|^2 \right) v \; ; \; v_{|t=0} = u_0 .$$

This is because $|u|^{2\sigma}$ and $|x|^{-\gamma} * |u|^2$ can be viewed as nonlinear potentials which do not see the phase term in (3.2), and because translation in the space variable commutes with the operation

$$u(t,x) \mapsto \left( |u|^{2\sigma}(t,x), \left( |x|^{-\gamma} * |u|^2 \right)(t,x) \right) .$$

When the potential is linear (a term $\mathtt{V}(x)u$ is added to the right hand side), this commutation property fails, and the role of the Stark potential is different (see the remark following Corollary 3.5).

The above result is formal. We give several corollaries where this property is used in cases where informations on $v^\varepsilon$ are available.

From now on, we take $\varepsilon = 1$ to simplify the presentation. The different results that we state are consequences of the above proposition and of known results for (1.3). These are not necessarily sharp results, but give a flavor of the properties of solutions to (1.1). We refer for instance to [4] and [9] for the results on (1.3).

**Corollary 3.2** (Local existence)**.** *Let $u_0 \in L^2(\mathbb{R}^n)$, $\lambda \in \mathbb{R}$ and $0 < \sigma < 2/n$. Then (1.1) has a unique solution*

$$u \in (C \cap L^\infty)(\mathbb{R}; L^2) \cap L^{2/\delta(2\sigma+2)}_{\mathrm{loc}}(\mathbb{R}; L^{2\sigma+2}) =: X_{2\sigma+2,\mathrm{loc}}(\mathbb{R}) ,$$

*which is actually in*

$$X_{\mathrm{loc}}(\mathbb{R}) := \left\{ u \in (C \cap L^\infty)(\mathbb{R}; L^2) \; ; \; u \in L^q_{\mathrm{loc}}(\mathbb{R}; L^r), \; \text{for all } 0 \leq \frac{2}{q} = \delta(r) < 1 \right\} ,$$

*and satisfies $\|u(t)\|_{L^2} = \|u_0\|_{L^2}$, $\forall t \in \mathbb{R}$.*

*Let $u_0 \in H^1(\mathbb{R}^n)$, $\lambda \in \mathbb{R}$ and $\sigma > 0$, with $\sigma < 2/(n-2)$ if $n \geq 3$. Then there exist $T_*(u_0), T^*(u_0) > 0$ such that (1.1) has a unique solution*

$$u \in X^1_{2\sigma+2,\mathrm{loc}}(]-T_*, T^*[) := \Big\{ u \in C(]-T_*, T^*[; L^2) \; ;$$

$$u, \nabla_x u \in C(]-T_*, T^*[; L^2) \cap L^{2/\delta(2\sigma+2)}_{\mathrm{loc}}(]-T_*, T^*[; L^{2\sigma+2}) \Big\},$$

*which is actually in $X^1_{\mathrm{loc}}(]-T_*, T^*[)$ (whose definition follows that of $X_{\mathrm{loc}}(\mathbb{R})$ and $X^1_{2\sigma+2,\mathrm{loc}}(]-T_*, T^*[)$). The following quantities are independent of $t \in ]-T_*, T^*[$,*

Mass: $\|u(t)\|_{L^2}$.

"Energy": $\dfrac{1}{2} \|(tE - i\nabla_x) u(t)\|^2_{L^2} + \dfrac{\lambda}{\sigma+1} \|u(t)\|^{2\sigma+2}_{L^{2\sigma+2}}$.

If $T^* < \infty$ (resp. $T_* < \infty$), then $\|\nabla_x u(t)\|_{L^2} \to \infty$ as $t \uparrow T^*$ (resp. as $t \downarrow -T_*$). If in addition $u_0 \in \Sigma$, then $u \in C(]-T_*, T^*[; \Sigma)$ and it satisfies (2.6).

*Remark.* Notice that using (3.2), we retrieve the definition of $J_E(t)$ from the usual Galilean operator $J(t) = x + it\nabla_x$, and Proposition 3.1 yields another interpretation of Lemmas 2.1 and 2.2.

*Remark.* The operator $tE - i\nabla_x$ is also an Heisenberg observable, and satisfies properties similar to those stated in Lemma 2.1 (see [3]).

It can be surprising to notice that the "natural" energy associated to (1.1),

$$(3.3) \qquad \frac{1}{2}\|\nabla_x u(t)\|_{L^2}^2 + \frac{\lambda}{\sigma+1}\|u(t)\|_{L^{2\sigma+2}}^{2\sigma+2} + \int V(x)|u(t,x)|^2 dx \, ,$$

which is formally conserved (see e.g. [18]), does not appear in the above result. We did not find a way to deduce this law from Proposition 3.1 and arguments assuming only that $u_0 \in H^1(\mathbb{R}^n)$ and $V|u_0|^2 \in L^1(\mathbb{R}^n)$. On the other hand, global existence results, as stated in the following corollary, do not seem easy to deduce from the conservation of (3.3).

Considering the two notions of energy introduced above, we see that the difference between these two quantities must be constant. This yields,

$$E \cdot \mathrm{Re} \int \overline{u}(t,x) J_E(t) u(t,x) dx = \mathrm{Const.}$$

Since (1.1) and (1.3) are related through (3.2), we infer that

$$E \cdot \mathrm{Re} \int \overline{v}(t,x) J(t) v(t,x) dx = \mathrm{Const.} \, ,$$

with $J(t) = x + it\nabla_x$. At this stage, $E$ is just a parameter and does not appear in (1.3). This implies

$$\mathrm{Re} \int \overline{v}(t,x) J(t) v(t,x) dx = \mathrm{Const.}$$

Notice that these identities were derived in [18], Claims 2 and 3.

**Corollary 3.3** (Global existence). *Let $u_0 \in H^1(\mathbb{R}^n)$, $\lambda \in \mathbb{R}$ and $\sigma > 0$, with $\sigma < 2/(n-2)$ if $n \geq 3$. We have $T^* = T_* = \infty$ in the following cases.*
- *The nonlinearity is repulsive, $\lambda \geq 0$.*
- *$\lambda < 0$ and $\sigma < 2/n$.*
- *$\lambda < 0$, $\sigma = 2/n$ and $\|u_0\|_{L^2} < \|Q\|_{L^2}$, where $Q$ is the unique spherically symmetric solution of (see [23], [15]),*

$$\begin{cases} -\frac{1}{2}\Delta Q + Q = -\lambda |Q|^{4/n} Q, \text{ in } \mathbb{R}^n, \\ \qquad\qquad\qquad Q > 0, \text{ in } \mathbb{R}^n. \end{cases}$$

- *$\lambda < 0$, $\sigma > 2/n$ and $\|u_0\|_{H^1}$ is sufficiently small.*

**Corollary 3.4** (Finite time blow up). *Let $u_0 \in \Sigma$, $\lambda < 0$ and $\sigma \geq 2/n$, with $\sigma < 2/(n-2)$ if $n \geq 3$. If*

$$\frac{1}{2}\|\nabla_x u_0\|_{L^2}^2 + \frac{\lambda}{\sigma+1}\|u_0\|_{L^{2\sigma+2}}^{2\sigma+2} < 0 \, ,$$

*then the solution $u$ to (1.1) blows up in finite time, in the future and in the past.*

*Remark.* It is well-known that finite time blow up occurs in other situations (see e.g. [16]). A remarkable fact is that from Proposition 3.1, the introduction of a linear potential (1.2) does not change the time of blow up, but only the location of this phenomenon (it has moved under the action of $V$ which can be interpreted as gravity). For quadratic potentials, not only the location, but also the time of blow up is altered (see [1], [2]).

It is well-known that $-(1/2)\Delta + V$ is essentially self-adjoint on $C_0^\infty(\mathbf{R^n})$ (see e.g. [14]). Let $H_S$ denote the closure of $-(1/2)\Delta + V$; we define a unitary group $U(t) = \exp(-itH_S)$. Then we state an immediate consequence of Proposition 3.1 on scattering theory (see e.g. [4], [5], [10], [11], [24], [26] and [5], [19] for the critical case); notice that the assumptions on $\sigma$ could be weakened for the existence of wave operators (first point). As mentioned above, we do not seek sharp results here.

**Corollary 3.5** (Scattering theory). *Let $u_0 \in \Sigma$, $\lambda > 0$, with $\sigma < 2/(n-2)$ if $n \geq 3$. Assume moreover that*
$$\sigma \geq \frac{2 - n + \sqrt{n^2 + 12n + 4}}{4n} \ .$$
- *For every $u_- \in \Sigma$, there exists a unique $u_0 \in \Sigma$ such that the maximal solution $u \in C(\mathbb{R}, \Sigma)$ to (1.1) satisfies*
$$\|U(-t)u(t) - u_-\|_\Sigma \underset{t \to -\infty}{\longrightarrow} 0.$$
- *For every $u_0 \in \Sigma$, there exists a unique $u_+ \in \Sigma$ such that the maximal solution $u \in C(\mathbb{R}, \Sigma)$ to (1.1) satisfies*
$$\|U(-t)u(t) - u_+\|_\Sigma \underset{t \to +\infty}{\longrightarrow} 0.$$

*Remark.* Proposition 3.1 shows that the critical indices for the power $\sigma$ in the scattering theory are the same as in the free case $E = 0$. For instance, in space dimension one, a short range theory is available if $\sigma > 1$, while a long range theory is needed if $\sigma \leq 1$ (see e.g. [9], [20], [12]; the expected critical index in space dimension $n$ is $\sigma_c = 1/n$). This is in sharp contrast with the linear case,
$$i\partial_t u + \frac{1}{2}\Delta u = E \cdot xu + |x|^{-\gamma}u \ , \ x \in \mathbb{R}^n \ .$$
If $E = 0$, and if $\gamma > 1$, then the potential $|x|^{-\gamma}$ is short range, while it is long range if $\gamma \leq 1$. On the other hand, if $E \neq 0$, then $\gamma = 1/2$ is critical. The different behaviors are due to the translation $x + (t^2/2)E$ (see for instance [21]).

**Acknowledgments**. The authors are grateful to Professor Jean Ginibre for pointing out reference [6]. The second author is grateful to Professor Tadayoshi Adachi for many helpful comments and encouragement.

(R. Carles) MAB, UMR CNRS 5466, Université Bordeaux 1, 351 cours de la Libération, 33 405 Talence cedex, France
*E-mail address*: carles@math.u-bordeaux.fr

(Y. Nakamura) Faculty of Engineering, Kumamoto University, 2-39-1, Kurokami, Kumamoto, 860-8555, Japan
*E-mail address*: hisa@math.sci.kumamoto.jp